\documentclass[12pt]{article}
\usepackage{amssymb}
\def\hybrid{\topmargin 0pt      \oddsidemargin 0pt
        \headheight 0pt \headsep 0pt
        \voffset=-0.5cm
        \hoffset=-0.25in
        \textwidth 6.75in
        \textheight 9.5in       % A4 paper
        \marginparwidth 0.0in
        \parskip 5pt plus 1pt   \jot = 1.5ex}
\catcode`\@=11
\def\marginnote#1{}

\newcount\hour
\newcount\minute
\newtoks\amorpm
\hour=\time\divide\hour by60 \minute=\time{\multiply\hour by60
\global\advance\minute by-\hour}
\edef\standardtime{{\ifnum\hour<12 \global\amorpm={am}%
        \else\global\amorpm={pm}\advance\hour by-12 \fi
        \ifnum\hour=0 \hour=12 \fi
        \number\hour:\ifnum\minute<10 0\fi\number\minute\the\amorpm}}
\edef\militarytime{\number\hour:\ifnum\minute<10 0\fi\number\minute}

\def\draftlabel#1{{\@bsphack\if@filesw {\let\thepage\relax
   \xdef\@gtempa{\write\@auxout{\string
      \newlabel{#1}{{\@currentlabel}{\thepage}}}}}\@gtempa
   \if@nobreak \ifvmode\nobreak\fi\fi\fi\@esphack}
        \gdef\@eqnlabel{#1}}
\def\@eqnlabel{}
\def\@vacuum{}
\def\draftmarginnote#1{\marginpar{\raggedright\scriptsize\tt#1}}
\def\draftlabel#1{{\@bsphack\if@filesw {\let\thepage\relax
   \xdef\@gtempa{\write\@auxout{\string
      \newlabel{#1}{{\@currentlabel}{\thepage}}}}}\@gtempa
   \if@nobreak \ifvmode\nobreak\fi\fi\fi\@esphack}
        \gdef\@eqnlabel{#1}}
\def\@eqnlabel{}
\def\@vacuum{}
\def\draftmarginnote#1{\marginpar{\raggedright\scriptsize\tt#1}}

\def\draft{\oddsidemargin -.5truein
        \def\@oddfoot{\sl preliminary draft \hfil
        \rm\thepage\hfil\sl\today\quad\militarytime}
        \let\@evenfoot\@oddfoot \overfullrule 3pt
        \let\label=\draftlabel
        \let\marginnote=\draftmarginnote
   \def\@eqnnum{(\theequation)\rlap{\kern\marginparsep\tt\@eqnlabel}%
\global\let\@eqnlabel\@vacuum}  }

\def\numberbysection{\@addtoreset{equation}{section}
        \def\theequation{\thesection.\arabic{equation}}}

\def\underline#1{\relax\ifmmode\@@underline#1\else
        $\@@underline{\hbox{#1}}$\relax\fi}

\def\titlepage{\@restonecolfalse\if@twocolumn\@restonecoltrue\onecolumn
     \else \newpage \fi \thispagestyle{empty}\c@page\z@
        \def\thefootnote{\fnsymbol{footnote}} }

\def\endtitlepage{\if@restonecol\twocolumn \else  \fi
        \def\thefootnote{\arabic{footnote}}
        \setcounter{footnote}{0}}  %\c@footnote\z@ }
\numberbysection \hybrid

\newcounter{mo}

\newcommand{\vf}{\varphi}
\newcommand{\al}{\alpha}
\newcommand{\be}{\beta}

\newcommand{\vth}{\vartheta}

\newcommand{\Mat}{ {\rm Mat}(N,\mathbb C) }

\newcommand{\mC}{\mathbb C}
\newcommand{\mZ}{\mathbb Z}

\def\beq{\begin{equation}}
\def\eq{\end{equation}}
\def\p{\partial}

\newtheorem{theor}{Theorem}%[section]

\def\res{\mathop{\hbox{Res}}\limits}

\newcommand{\ox}{\otimes}

\begin{document}

\setcounter{page}{1}

\date{}
\date{}
\vspace{50mm}

\begin{flushright}
 ITEP-TH-29/15\\
\end{flushright}
\vspace{5mm}

\begin{center}
\vspace{0mm}
{\LARGE{Associative Yang-Baxter equation for }}
 \\ \vspace{4mm}
 {\LARGE{quantum (semi-)dynamical R-matrices}}
\\
\vspace{12mm} {\large \ \ \ \ {Ivan Sechin}\,$^{\natural\, \sharp}$
\ \ \ \
 \ \ \ {Andrei Zotov}\,$^{\diamondsuit\, \sharp\, \natural}$ }\\
 \vspace{8mm}
% \vspace{3mm}

 \vspace{2mm} $^\natural$ -- {\small{\sf MIPT, Inststitutskii per.  9, Dolgoprudny,
 Moscow region, 141700, Russia}}\\
 \vspace{2mm} $^\sharp$ -- {\small{\sf %Institute of Theoretical and Experimental Physics, 117218,  Moscow, Russia
 ITEP, B. Cheremushkinskaya str. 25,  Moscow, 117218, Russia}}\\
\vspace{2mm} $^\diamondsuit$ -- {\small{\sf Steklov Mathematical
Institute  RAS, Gubkina str. 8, Moscow, 119991,  Russia}}
\end{center}

%\vspace{2mm}
\begin{center}\footnotesize{{\rm E-mails:}{\rm\ \
 shnbuz@gmail.com,\ \ zotov@mi.ras.ru}}\end{center}
%
%\vspace{0mm}
%
 \begin{abstract}
In this paper we propose versions of the associative Yang-Baxter
equation and higher order $R$-matrix identities which can be applied
to quantum dynamical $R$-matrices. As is known quantum non-dynamical
$R$-matrices of Baxter-Belavin type satisfy this equation. Together
with unitarity condition and skew-symmetry it provides the quantum
Yang-Baxter equation and a set of identities useful for different
applications in integrable systems. The dynamical $R$-matrices
satisfy the Gervais-Neveu-Felder (or dynamical Yang-Baxter)
equation.  Relation between the dynamical and non-dynamical cases is
described by the IRF-Vertex transformation. An alternative approach
to quantum (semi-)dynamical $R$-matrices and related quantum
algebras was suggested by Arutyunov, Chekhov and Frolov (ACF) in
their study of the quantum Ruijsenaars-Schneider model. The purpose
of this paper is twofold. First, we prove that the ACF elliptic
$R$-matrix satisfies the associative Yang-Baxter equation with
shifted spectral parameters. Second, we directly prove a simple
relation of the IRF-Vertex type between the Baxter-Belavin and the
ACF elliptic $R$-matrices predicted previously by Avan and Rollet.
It provides the higher order $R$-matrix identities and an
explanation of the obtained equations through those for
non-dynamical $R$-matrices. As a by-product we also get an
interpretation of the intertwining transformation as matrix
extension of scalar theta function likewise $R$-matrix is
interpreted as matrix extension of the Kronecker function. Relations
to the Gervais-Neveu-Felder equation and identities for the Felder's
elliptic $R$-matrix are also discussed.
 \end{abstract}

\newpage

{\small{

\tableofcontents

}}

%\small{

\section{Introduction and summary}
\setcounter{equation}{0}
We start with a brief review of the Yang-Baxter structures under
consideration, and then give a summary of the paper.
\subsection{Brief review}
\noindent {\bf Quantum non-dynamical $R$-matrices} of ${\rm
GL}(N,\mC)$ type in fundamental representation are elements of
$\Mat^{\otimes 2}$ satisfying the quantum Yang-Baxter equation
\cite{Yang,Baxter,Sklyanin}
  \beq\label{a01}
  \begin{array}{c}
  \displaystyle{
 R^\hbar_{12}(z_1,z_2)R^\hbar_{13}(z_1,z_3)R^\hbar_{23}(z_2,z_3)
 =R^\hbar_{23}(z_2,z_3)R^\hbar_{13}(z_1,z_3)R^\hbar_{12}(z_1,z_2)
 }
 \end{array}
 \eq
together with the unitarity condition which we write in some special
$R$-matrix normalization\footnote{The normalization (\ref{a02})
implies in fact that we deal with a special class of $R$-matrices
which includes Baxter-Belavin's elliptic one and some of its
trigonometric and rational degenerations.}:
  \beq\label{a02}
  \begin{array}{c}
  \displaystyle{
R^\hbar_{12}(z_1,z_2)
R^\hbar_{21}(z_2,z_1)=\phi(\hbar,z_1-z_2)\phi(\hbar,z_2-z_1)\,\,1\otimes
1= (\wp(\hbar)-\wp(z_1-z_2))\,\,1\otimes 1\,,
 }
 \end{array}
 \eq
where $\phi(\hbar,z)$ is the Kronecker function (\ref{a907}) and
$\wp(z)$ is the Weierstrass $\wp$-function (\ref{a908}),
(\ref{a912}). The parameter $\hbar$ is called the Planck constant,
and $z_1,z_2$ -- spectral parameters. The Baxter-Belavin
\cite{Baxter,Belavin} elliptic solution of (\ref{a01}), (\ref{a02})
is of the form:
 \beq\label{a03}
 \begin{array}{c}
  \displaystyle{
 R^{\hbox{\tiny{B}}}_{12}(\hbar,z_1,z_2)= R^{\hbox{\tiny{B}}}_{12}(\hbar,z_1-z_2)=
 \frac{1}{N}\sum\limits_{a\in\,{\mathbb Z}_N\times{\mathbb Z}_N} \vf_a^\hbar(z_1-z_2)\, T_a\otimes
 T_{-a}\,,
 }
 \end{array}
 \eq
where $\{T_a\}$ is a special basis (\ref{a904}) in $\Mat$, and
$\{\vf_a^\hbar(z)\}$ is a set of related functions (\ref{a910}).

\vskip3mm

\noindent {\bf Quantum dynamical $R$-matrices: Gervais-Neveu-Felder
equation.} The dynamical $R$-matrices depend on additional
parameters $u_1,...,u_N$. In the classical Hamiltonian mechanics of
integrable many-body systems they are the coordinates of particles,
while in quantum case these are the parameters entering the
Boltzmann weights in IRF statistical models \cite{Jimbo,Jimbo2}.

The quantum dynamical $R$-matrices are described by the
Gervais-Neveu-Felder equation \cite{Felder}:
  \beq\label{a04}
  \begin{array}{c}
  \displaystyle{
 R^\hbar_{12}(z_1,z_2|\,u)R^\hbar_{13}(z_1,z_3|\,u+\hbar^{(2)})R^\hbar_{23}(z_2,z_3|\,u)=\hspace{40mm}
}
\\ \ \\
  \displaystyle{
\hspace{40mm}
=R^\hbar_{23}(z_2,z_3|\,u+\hbar^{(1)})R^\hbar_{13}(z_1,z_3|\,u)R^\hbar_{12}(z_1,z_2|\,u+\hbar^{(3)})\,,
 }
 \end{array}
 \eq
where the shifts of dynamical arguments $u$ are performed as
follows:
  \beq\label{a041}
  \begin{array}{c}
  \displaystyle{
R^\hbar_{12}(z_1,z_2|\,u+\hbar^{(3)})=P_3^\hbar\,
R^\hbar_{12}(z_1,z_2|\,u)\, P_3^{-\hbar} \,,\quad
P_3^\hbar=\sum\limits_{k=1}^N 1\otimes 1\otimes E_{kk}
\exp(\hbar\frac{\p}{\p u_k})
 }
 \end{array}
 \eq
with notation $\{E_{ij}\}$ for the standard basis in $\Mat$:
$(E_{ij})_{kl}=\delta_{ik}\delta_{jl}$. The weight zero condition
implies that $[R^\hbar_{12}(z_1,z_2|\,u),P_1^\hbar P_2^\hbar]=0$.

The elliptic solution of (\ref{a04}), (\ref{a02}) is given by the
Felder's $R$-matrix \cite{Felder2}:
 \beq\label{a05}
 \begin{array}{c}
  \displaystyle{
 R^{\hbox{\tiny{F}}}_{12}(\hbar,z_1,z_2|\,u)=R^{\hbox{\tiny{F}}}_{12}(\hbar,z_1-z_2|\,u)=
 }
\\ \ \\
  \displaystyle{
 =\sum\limits_{i\neq j}
 E_{ii}\otimes E_{jj}\, \phi(\hbar,u_{ij})+\sum\limits_{i\neq j}
 E_{ij}\otimes E_{ji}\, \phi(z_1-z_2,-u_{ij})+\phi(\hbar,z_1-z_2)\sum\limits_{i}
 E_{ii}\otimes E_{ii}\,,
 }
 \end{array}
 \eq
where $u_{ij}=u_i-u_j$ and $\phi(\eta,z)$ is the Kronecker function
(\ref{a907}).

\vskip3mm

\noindent {\bf IRF-Vertex (or Vertex-Face) correspondence} provides
explicit relation between dynamical and non-dynamical $R$-matrices
\cite{Baxter2,Jimbo,Hasegawa}. Its applications to classical
integrable systems, 1+1 models and monodromy preserving equations
can be found in \cite{LOZ,LOZ8}. Consider the following matrix
$g\in\Mat$:
 \beq\label{a06}
 \begin{array}{c}
  \displaystyle{
g_{ij}(z,u)=
 \vth\left[  \begin{array}{c}
 \frac12-\frac{i}{N} \\ \frac12
 \end{array} \right] (z+Nu_j-\sum\limits_{m=1}^N
 u_m|\,N\tau)\frac{1}{\prod\limits_{k\neq j}\vth(u_k-u_j)}\,,
 }
 \end{array}
 \eq
where theta-functions with characteristics are defined in
(\ref{a901}). The IRF-Vertex relation between the quantum
$R$-matrices (\ref{a03}) and (\ref{a05}) has the form:
 \beq\label{a07}
 \begin{array}{c}
  \displaystyle{
g_2(z_2,u)\,
g_1(z_1,u-\hbar^{(2)})\,R^{\hbox{\tiny{F}}}_{12}(\hbar,z_1-z_2|\,u)=R^{\hbox{\tiny{B}}}_{12}(\hbar,z_1-z_2)\,
g_1(z_1,u)\, g_2(z_2,u-\hbar^{(1)})\,.
 }
 \end{array}
 \eq
The trigonometric and rational analogues of (\ref{a03}) and
(\ref{a06}) can be found in \cite{Zabr} and \cite{LOZ8}.

\vskip3mm

\noindent {\bf Quantum dynamical $R$-matrices:
Arutyunov-Chekhov-Frolov (ACF) approach.} An alternative to
(\ref{a04}) quantization of dynamical $r$-matrix structure
%(related to the Ruijsenaars-Schneider model)
was suggested in \cite{Arut} and then studied in \cite{Avan1,Avan2},
where it was called semi-dynamical Yang-Baxter equation:
  \beq\label{a08}
  \begin{array}{c}
  \displaystyle{
 R^\hbar_{12}(z_1,z_2|\,u)R^\hbar_{13}(z_1-\hbar,z_3-\hbar|\,u)R^\hbar_{23}(z_2,z_3|\,u)=\hspace{50mm}
}
\\ \ \\
  \displaystyle{
\hspace{40mm}
=R^\hbar_{23}(z_2-\hbar,z_3-\hbar|\,u)R^\hbar_{13}(z_1,z_3|\,u)R^\hbar_{12}(z_1-\hbar,z_2-\hbar|\,u)\,.
 }
 \end{array}
 \eq
Notice that in (\ref{a08}) there are no shifts (\ref{a041}) of the
dynamical parameters but there are shifts of the spectral parameters
instead. The elliptic $R$-matrix satisfying the unitarity condition
(\ref{a02}) and the Yang-Baxter equation (\ref{a08}) was found in
\cite{Arut}. It is of the form:\footnote{We use different sign for
the dynamical parameters $u$, and $R$-matrix normalization is chosen
as in (\ref{a02}).}
 \beq\label{a09}
 \begin{array}{c}
  \displaystyle{
 R^{\hbox{\tiny{ACF}}}_{12}(\hbar,z_1,z_2|\,u)=\sum\limits_{i\neq j}
 E_{ii}\otimes E_{jj}\, \phi(\hbar,-u_{ij})+\sum\limits_{i\neq j}
 E_{ij}\otimes E_{ji}\, \phi(z_1-z_2,-u_{ij})-
 }
\\ \ \\
  \displaystyle{
-\sum\limits_{i\neq j}
 E_{ij}\otimes E_{jj}\, \phi(z_1+\hbar,-u_{ij})+\sum\limits_{i\neq j}
 E_{jj}\otimes E_{ij}\, \phi(z_2,-u_{ij})+
 }
\\ \ \\
  \displaystyle{
+(E_1(\hbar)+E_1(z_1-z_2)+E_1(z_2)-E_1(z_1+\hbar))\sum\limits_{i}
 E_{ii}\otimes E_{ii}\,,
 }
 \end{array}
 \eq
where $u_{ij}=u_i-u_j$ and $E_1(z)=\vth'(z)/\vth(z)$ is the first
Eisenstein function (\ref{a908}).

The relation between $R^{\hbox{\tiny{ACF}}}$ and
$R^{\hbox{\tiny{F}}}$ was obtained
  \beq\label{a10}
  \begin{array}{c}
  \displaystyle{
 R^{\hbox{\tiny{ACF}}}_{12}(\hbar,z_1,z_2|\,u)={\bar
 R}_{12}(\hbar,z_1|\,u-\hbar^{(2)})\,
 R^{\hbox{\tiny{F}}}_{12}(\hbar,z_1-z_2|\,u)\, {\bar
 R}_{21}^{-1}(\hbar,z_2|\,u-\hbar^{(1)})
 }
 \end{array}
 \eq
 or
  \beq\label{a101}
  \begin{array}{c}
  \displaystyle{
 R^{\hbox{\tiny{ACF}}}_{12}(\hbar,z_1,z_2|\,u)={\bar
 R}_{21}(\hbar,z_2|\,u)\,
 R^{\hbox{\tiny{F}}}_{12}(\hbar,z_1-z_2|\,u)\,
 {\bar
 R}_{12}^{-1}(\hbar,z_1|\,u)
 }
 \end{array}
 \eq
 in terms of
explicitly given twist matrix\footnote{The shifts $\hbar^{(1)}$,
$\hbar^{(2)}$ in (\ref{a10}) are defined as in (\ref{a041}).}
  \beq\label{a11}
  \begin{array}{l}
  \displaystyle{
 \frac{\vth'(0)}{\vth(\hbar)}\,{\bar R}_{12}(\hbar,z|\,u)=
 }
\\ \ \\
  \displaystyle{
 =\sum\limits_{i\neq j}
 E_{ii}\otimes E_{jj}\, \phi(\hbar,-u_{ij})-\sum\limits_{i\neq j}
 E_{ij}\otimes E_{jj}\, \phi(z+\hbar,-u_{ij})-\phi(z+\hbar,-\hbar)\sum\limits_{i}
 E_{ii}\otimes E_{ii}
 }
 \end{array}
 \eq
 and its inverse
  \beq\label{a12}
  \begin{array}{l}
  \displaystyle{
 \frac{\vth'(0)}{\vth(\hbar)}\,{\bar R}_{12}^{-1}(\hbar,z|\,u)=
\sum\limits_{i,j}
 E_{ii}\otimes E_{jj}\, \phi(\hbar,u_{ij}-\hbar)-\sum\limits_{i, j}
 E_{ij}\otimes E_{jj}\, \phi(z,\hbar-u_{ij})\,.
 }
 \end{array}
 \eq
The origin of the ACF type Yang-Baxter equation (\ref{a08}), i.e.
its relation to the Yang-Baxter equation (\ref{a01}) or the
Gervais-Neveu-Felder equation (\ref{a04}) was described in
\cite{Avan1,Avan2}\footnote{A certain relation between (\ref{a08})
and (\ref{a01}) was observed indirectly in the original paper
\cite{Arut}: it was shown that any representation of quantum algebra
underlying $R^{\hbox{\tiny{ACF}}}$:
$R^{\hbox{\tiny{ACF}}}_{12}(z,w)L_1(z){\bar
R}_{21}(w)L_2(w)=L_2(w){\bar
R}_{12}(z)L_1(z)R^{\hbox{\tiny{F}}}_{12}(z,w)$ turns into a
representation of the exchange relations
$R^{\hbox{\tiny{B}}}_{12}(z-w){\hat L}_1(z){\hat L}_2(w)={\hat
L}_2(w){\hat L}_1(z)R^{\hbox{\tiny{B}}}_{12}(z-w)$ via the
IRF-Vertex transformation (\ref{a07}).}. We will give a direct proof
of relation between $R$-matrices $R^{\hbox{\tiny{ACF}}}$ and
$R^{\hbox{\tiny{B}}}$ (and therefore $R^{\hbox{\tiny{F}}}$) and
corresponding Yang-Baxter equations as well as higher order
$R$-matrix identities (see Theorem \ref{theor2} below).

\vskip3mm

\noindent {\bf Associative Yang-Baxter equation}
%for non-dynamical quantum $R$-matrices. The equation
was originally introduced in \cite{Aguiar} for constant $R$-matrices
and then generalized by Polishchuk \cite{Pol} to the form:
  \beq\label{a15}
  \begin{array}{c}
  \displaystyle{
 R^\hbar_{12}
 R^{\eta}_{23}=R^{\eta}_{13}R_{12}^{\hbar-\eta}+R^{\eta-\hbar}_{23}R^\hbar_{13}\,,\
 \ R^\hbar_{ab}=R^\hbar_{ab}(z_a\!-\!z_b)\,.
 }
 \end{array}
 \eq
%
% triple Massey products and homological mirror symmetry on elliptic curves.
%
It was shown in \cite{Pol} that the Baxter-Belavin $R$-matrix
written in Richey-Tracy \cite{RicheyT} form satisfies (\ref{a15}).
The unitarity condition (\ref{a02}) was not required. The $R$-matrix
(\ref{a03}) was rather considered as deformation of the classical
one. In fact, the quantum $R$-matrix (\ref{a03}) satisfies also the
skew-symmetry property
  \beq\label{a16}
  \begin{array}{c}
  \displaystyle{
 R^\hbar_{12}(z_1-z_2)=-R_{21}^{-\hbar}(z_2-z_1)
 }
 \end{array}
 \eq
 which can be viewed as the classical analogue of the unitarity
condition (\ref{a02}).\footnote{The interpretation of (\ref{a03}) as
the classical $r$-matrix is also discussed in \cite{LOZ11} briefly.}
The relation of (\ref{a15}) to the quantum Yang-Baxter equation
(\ref{a01}) was studied separately \cite{Pol3}.

The most natural\footnote{It is natural from the viewpoint of the
underlying Fay identity. See (\ref{a21}) below.} and simple
dynamical solution of (\ref{a15}) was proposed by Burban and Henrich
\cite{Burban2} (see also \cite{Rubt}):
 \beq\label{a17}
 \begin{array}{c}
  \displaystyle{
 R^{\hbox{\tiny{BH}}}_{12}(\hbar,z_1,z_2|\,u)=R^{\hbox{\tiny{BH}}}_{12}(\hbar,z_1-z_2|\,u)=
\sum\limits_{i, j}
 E_{ij}\otimes E_{ji}\, \phi(z_1-z_2,\hbar-u_{ij})\,.
 }
 \end{array}
 \eq
 It is skew-symmetric (\ref{a16}) but the unitarity condition is not
 valid:
 \beq\label{a18}
 \begin{array}{c}
  \displaystyle{
R^{\hbox{\tiny{BH}}}_{12}(\hbar,z_1-z_2|\,u)R^{\hbox{\tiny{BH}}}_{21}(\hbar,z_2-z_1|\,u)=
\sum\limits_{i, j}
 E_{ii}\otimes E_{jj}\, (\wp(\hbar-u_{ij})-\wp(z_1-z_2))
 }
 \end{array}
 \eq
 as well as the quantum Yang-Baxter equation (\ref{a01}).
Therefore, it behaves more like a classical
$r$-matrix.\footnote{However it may be useful for applications to
integrable systems. We discuss it in our future papers, where
$R$-matrix valued Lax pairs and KZB equations related to quantum
dynamical (and/or semi-dynamical) $R$-matrices will be studied as
well as possible generalization of the IRF-Vertex transformations
for these structures.}

Later \cite{LOZ9,LOZ10} the equation (\ref{a15}) found applications
in integrable systems, the KZB and Painlev\'e equations. In
particular, it was mentioned in \cite{LOZ9,LOZ11} that a unitary
(\ref{a02}) and skew-symmetric (\ref{a16}) solution of the
associative Yang-Baxter equation (\ref{a15}) (in particular, the
Baxter-Belavin (\ref{a03}) one) satisfies also the following cubic
identity\footnote{The cubic relation (\ref{a19}) follows from
(\ref{a15}) and (\ref{a02}), (\ref{a16}) but the inverse statement
is unknown.}:
  \beq\label{a19}
  \begin{array}{c}
  \displaystyle{
 R^\eta_{12} R^\hbar_{13} R^\eta_{23}-R^\hbar_{23}
 R^\eta_{13} R^\hbar_{12}=R_{13}^{\hbar+\eta}\,\left( \wp(\eta)-\wp(\hbar)
 \right)\,,
 }
 \end{array}
 \eq
where $R^\hbar_{ab}=R^{\hbox{\tiny{B}}}_{ab}(\hbar,z_a,z_b)$. For
$\eta=\hbar$ the latter equation provides the Yang-Baxter one
(\ref{a01}) while for $\eta=-\hbar$ it leads to
  \beq\label{a20}
  \begin{array}{c}
  \displaystyle{
 R^\hbar_{12} R^\hbar_{23} R^\hbar_{31}+R^\hbar_{13}
 R^\hbar_{32} R^\hbar_{21}=-\wp'(\hbar)\,\,1\otimes 1\otimes 1\,.
 }
 \end{array}
 \eq
Higher order analogues of (\ref{a20}) can be found in \cite{Z}. They
are discussed below.

\subsection{Summary}

The purpose of paper is to find solutions of the associative
Yang-Baxter equation (\ref{a15}) and to prove identities of type
(\ref{a19}), (\ref{a20}) for quantum (semi-)dynamical $R$-matrices.
Before we proceed further let us mention that  (\ref{a15}) can be
considered as matrix generalization of the Fay identity
(\ref{a909}), which for $z=z_1-z_2$ and $w=z_2-z_3$ takes the form
  \beq\label{a21}
  \begin{array}{c}
  \displaystyle{
\phi(\hbar,z_{12})\phi(\eta,z_{23})=\phi(\eta,z_{13})\phi(\hbar-\eta,z_{12})+\phi(\eta-\hbar,z_{23})\phi(\hbar,z_{13})\,,
\quad z_{ab}=z_a-z_b\,.
 }
 \end{array}
 \eq
Indeed, the Baxter-Belavin $R$-matrix in scalar case (for $N=1$) is
exactly the Kronecker function $\phi(\hbar,z_1-z_2)$. This analogy,
in fact, underlies the results of papers \cite{Pol} and
\cite{LOZ9,LOZ10,Z}. In this sense the unitarity condition
(\ref{a02}) is analogue of (\ref{a912}).

Similarly to the Baxter-Belavin case the Felder's $R$-matrix
(\ref{a05}) is also unitary and skew-symmetric. Moreover, for $N=1$
it equals to the Kronecker function $\phi(\hbar,z_1-z_2)$. However,
we have not found a quadratic equation of type (\ref{a15}) for
$R^{\hbox{\tiny{F}}}$. Equations for $R^{\hbox{\tiny{F}}}$ follows
from those for for $R^{\hbox{\tiny{B}}}$ via the IRF-Vertex
transformation (\ref{a07}) but the twist matrix (\ref{a06}) is not
cancelled out from the final answers. We discuss it in Section
\ref{sect4}. At the same time it appears that the quadratic equation
of type (\ref{a15}) holds true for $R^{\hbox{\tiny{ACF}}}$.

\begin{theor}\label{theor1}
 The ACF $R$-matrix (\ref{a09}) satisfies the following modification
 of the associative Yang-Baxter equation (\ref{a15}):
  \beq\label{a23}
  \begin{array}{c}
  \displaystyle{
 R^\hbar_{12}(z_1+\eta,z_2+\eta)
 R^{\eta}_{23}(z_2+\hbar,z_3+\hbar)=
 }
 \\ \ \\
  \displaystyle{
 R^{\eta}_{13}(z_1+\hbar,z_3+\hbar)R_{12}^{\hbar-\eta}(z_1+\eta,z_2+\eta)
 +R^{\eta-\hbar}_{23}(z_2+\hbar,z_3+\hbar)R^\hbar_{13}(z_1+\eta,z_3+\eta)\,,
 }
 \end{array}
 \eq
where $R^\hbar_{ab}(z,w)=R^{\hbox{\tiny{ACF}}}_{ab}(\hbar,z,w|\,u)$,
and the cubic identity
  \beq\label{a24}
  \begin{array}{c}
  \displaystyle{
 R^\eta_{12}(z_1,z_2)\, R^\hbar_{13}(z_1\!-\!\hbar,z_3\!-\!\hbar)\,
 R^\eta_{23}(z_2,z_3)-
 }
 \\ \ \\
  \displaystyle{
-R^\hbar_{23}(z_2\!-\!\hbar,z_3\!-\!\hbar)\, R^\eta_{13}(z_1,z_3)\,
R^\hbar_{12}(z_1\!-\!\hbar,z_2\!-\!\hbar)
=R_{13}^{\hbar+\eta}(z_1\!-\!\hbar,z_3\!-\!\hbar)\,\left(
\wp(\eta)-\wp(\hbar)
 \right)\,.
 }
 \end{array}
 \eq
\end{theor}
See Section \ref{sect2} for the proof. As a conclusion of this
theorem we also obtain the Yang-Baxter equation (\ref{a08})
($\eta=\hbar$ in (\ref{a24})) and the unchanged identity (\ref{a19})
($\eta=-\hbar$ in (\ref{a24})).

The results of Theorem \ref{theor1} are valid in trigonometric and
rational cases as well. The functions (\ref{a907}), (\ref{a908})
entering the ACF $R$-matrix are given for these cases.

The equations (\ref{a23}) and  (\ref{a24}) can be derived from their
non-dynamical analogues (\ref{a15}) and (\ref{a19}) using the
IRF-Vertex like relation between $R^{\hbox{\tiny{ACF}}}$ (\ref{a09})
and $R^{\hbox{\tiny{B}}}$ (\ref{a03}) predicted in \cite{Avan1}:
 \begin{theor}\label{theor2}
The ACF $R$-matrix (\ref{a09}) and the twist matrix (\ref{a11}) are
expressed in terms of IRF-Vertex transformation matrix (\ref{a06})
and the Baxter-Belavin $R$-matrix (\ref{a03}) as follows:
  \beq\label{a25}
  \begin{array}{c}
  \displaystyle{
 R^{\hbox{\tiny{B}}}_{12}(\hbar,z_1-z_2)=g_1(z_1+\hbar,u)\,
 g_2(z_2,u)\,
 R^{\hbox{\tiny{ACF}}}_{12}(\hbar,z_1,z_2|\,u)\, g_2^{-1}(z_2+\hbar,u)\,
 g_1^{-1}(z_1,u)
 }
 \end{array}
 \eq
 and
   \beq\label{a26}
  \begin{array}{c}
  \displaystyle{
 {\bar R}_{12}(\hbar,z|\,u)=g_1^{-1}(z+\hbar,u+\hbar^{(2)})\,
 g_1(z,u)\,.
 }
 \end{array}
 \eq
 \end{theor}
See the proof in Section \ref{sect3}. As we will see it follows from
(\ref{a25}) that the ACF $R$-matrix satisfies also $n$-th order
identities proved for $R^{\hbox{\tiny{B}}}$ in \cite{Z}:
  \beq\label{a27}
  \begin{array}{c}
  \displaystyle{
 \sum\limits_{\hbox{\tiny{$
\begin{array}{c}
1\leq i_1 ... i_{n\!-\!1}\leq n
\\
i_c\neq a;\ i_b\neq i_c
\end{array}
$}}}
 R^\hbar_{ai_1} R^\hbar_{i_1 i_2}\,...\, R^\hbar_{i_{n-2}i_{n-1}} R^\hbar_{i_{n-1}a}=
 \underbrace{1\otimes...\otimes 1}_{\hbox{\small{\em n
times}}}\,(-1)^n\left.\frac{d^{(n-2)}}{d\eta^{(n-2)}}\,\wp(\eta)\right|_{\eta=\hbar}\,.
 }
 \end{array}
 \eq
 where $R^\hbar_{ij}=R^{\hbox{\tiny{ACF}}}_{ij}(\hbar,z_i,z_j|\,u)$, $a$ is a fixed index $1\leq a \leq n$ and $n\geq 3$. For
 $n=3$ it is (\ref{a20}).

As already mentioned the associative Yang-Baxter equation
(\ref{a15}) in non-dynamical case is a matrix analogue of the Fay
identity (\ref{a21}) for the Kronecker function. At the same time
the scalar ($N=1$) ACF $R$-matrix is not the Kronecker function. It
is equal to
  \beq\label{a28}
  \begin{array}{c}
  \displaystyle{
\varrho^\hbar(z_1,z_2)=E_1(\hbar)+E_1(z_1-z_2)+E_1(z_2)-E_1(z_1+\hbar)\stackrel{(\ref{a911})}{=}
 \phi(\hbar,z_1-z_2)\,\frac{\phi(\hbar,z_2)}{\phi(\hbar,z_1)}\,.
 }
 \end{array}
 \eq
Nevertheless this function satisfies (\ref{a23}) due to the Fay
identity (\ref{a21}) because the products of additional multiples
(${\phi(\hbar,z_i)}/{\phi(\hbar,z_j)}$) are equal for each term of
(\ref{a23}). See (\ref{a99}).
%
%  $$
%  \begin{array}{c}
%  \displaystyle{
% \frac{\phi(\hbar,z_2+\eta)}{\phi(\hbar,z_1+\eta)}\,\frac{\phi(\eta,z_3+\hbar)}{\phi(\eta,z_2+\hbar)}=
% \frac{\phi(\eta,z_3+\hbar)}{\phi(\eta,z_1+\hbar)}\,\frac{\phi(\hbar-\eta,z_2+\eta)}{\phi(\hbar-\eta,z_1+\eta)}=
% \frac{\phi(\eta-\hbar,z_3+\hbar)}{\phi(\eta-\hbar,z_2+\hbar)}\,
% \frac{\phi(\hbar,z_3+\eta)}{\phi(\hbar,z_1+\eta)}\,.
% }
% \end{array}
% $$

It is also notable that (\ref{a25}) leads to interpretation of the
intertwining matrix (\ref{a06}) as a matrix analogue of theta
function similarly to interpretation of $R$-matrix (\ref{a03}) as a
matrix analogue of the Kronecker function. See (\ref{a64}).
Finally, in Section \ref{sect4} we discuss identities for the
Felder's $R$-matrix arising from the IRF-Vertex relations.

\section{Associative YB equation for ACF $R$-matrix}\label{sect2}
\setcounter{equation}{0}

The proof of (\ref{a23}) is given in the Appendix. Let us derive the
cubic identity (\ref{a24}) likewise it was made in \cite{LOZ11} for
derivation of (\ref{a19}) using (\ref{a15}), (\ref{a02}) and
(\ref{a16}). For this purpose we also need an analogue of the
skew-symmetry property (\ref{a16}) for $R^{\hbox{\tiny{ACF}}}$. It
is given in \cite{Arut}:
  \beq\label{a41}
  \begin{array}{l}
  \displaystyle{
R^\hbar_{12}(z_1,z_2)=-R^{-\hbar}_{21}(z_2+\hbar,z_1+\hbar)\,,
 }
 \end{array}
 \eq
 where
 $R^\hbar_{12}(z_1,z_2)=R^{\hbox{\tiny{ACF}}}_{12}(\hbar,z_1,z_2|\,u)$.

\noindent \vskip1mm\underline{\em{Proof of cubic identity
(\ref{a24})}}:\vskip1mm

Multiplying equation (\ref{a23}) by
$R^{\hbar-\eta}_{23}(z_2+\eta,z_3+\eta)$ form the left and using
(\ref{a02}), (\ref{a41}) we obtain
  \beq\label{a42}
  \begin{array}{c}
  \displaystyle{
R^{\hbar-\eta}_{23}(z_2+\eta,z_3+\eta)
 R^\hbar_{12}(z_1+\eta,z_2+\eta)
 R^{\eta}_{23}(z_2+\hbar,z_3+\hbar)=
 }
 \\ \ \\
  \displaystyle{
=R^{\hbar-\eta}_{23}(z_2+\eta,z_3+\eta)
R^{\eta}_{13}(z_1+\hbar,z_3+\hbar)R_{12}^{\hbar-\eta}(z_1+\eta,z_2+\eta)-
 }
 \\ \ \\
  \displaystyle{
 -(\wp(\hbar-\eta)-\wp(z_2-z_3))\, R^\hbar_{13}(z_1+\eta,z_3+\eta)\,.
 }
 \end{array}
 \eq
Consider (\ref{a23}) with interchanged indices 2 and 3:
  \beq\label{a43}
  \begin{array}{c}
  \displaystyle{
 R^\hbar_{13}(z_1+\eta,z_3+\eta)
 R^{\eta}_{32}(z_3+\hbar,z_2+\hbar)=
 }
 \\ \ \\
  \displaystyle{
 R^{\eta}_{12}(z_1+\hbar,z_2+\hbar)R_{13}^{\hbar-\eta}(z_1+\eta,z_3+\eta)
 +R^{\eta-\hbar}_{32}(z_3+\hbar,z_2+\hbar)R^\hbar_{12}(z_1+\eta,z_2+\eta)\,.
 }
 \end{array}
 \eq
Multiplying it by $R^{\eta}_{23}(z_2+\hbar,z_3+\hbar)$ from the
right and using (\ref{a02}), (\ref{a41}) we get:
  \beq\label{a44}
  \begin{array}{c}
  \displaystyle{
 R^\hbar_{13}(z_1+\eta,z_3+\eta)\,
 (\wp(\eta)-\wp(z_2-z_3))=
 }
 \\ \ \\
  \displaystyle{
 =R^{\eta}_{12}(z_1+\hbar,z_2+\hbar)R_{13}^{\hbar-\eta}(z_1+\eta,z_3+\eta)R^{\eta}_{23}(z_2+\hbar,z_3+\hbar)+
 }
 \\ \ \\
  \displaystyle{
 +R^{\eta-\hbar}_{32}(z_3+\hbar,z_2+\hbar)R^\hbar_{12}(z_1+\eta,z_2+\eta)R^{\eta}_{23}(z_2+\hbar,z_3+\hbar)\,.
 }
 \end{array}
 \eq
Subtracting (\ref{a44}) from (\ref{a42}) yields:
  \beq\label{a45}
  \begin{array}{c}
  \displaystyle{
 R^{\eta}_{12}(z_1+\hbar,z_2+\hbar)R_{13}^{\hbar-\eta}(z_1+\eta,z_3+\eta)R^{\eta}_{23}(z_2+\hbar,z_3+\hbar)-
 }
 \\ \ \\
  \displaystyle{
 -R^{\hbar-\eta}_{23}(z_2+\eta,z_3+\eta)
R^{\eta}_{13}(z_1+\hbar,z_3+\hbar)R_{12}^{\hbar-\eta}(z_1+\eta,z_2+\eta)=
 }
 \\ \ \\
  \displaystyle{
  (\wp(\eta)-\wp(\hbar-\eta))\,R^\hbar_{13}(z_1+\eta,z_3+\eta)\,.

 }
 \end{array}
 \eq
Redefinition $\hbar:=\hbar+\eta$ and
$z_{1,2,3}:=z_{1,2,3}-\hbar-\eta$ gives (\ref{a24}). $\blacksquare$

The special case $\eta=\hbar$ for (\ref{a24}) obviously reproduces
the Yang-Baxter equation (\ref{a08}). At the same time the case
$\eta=-\hbar$ yields (\ref{a20})
  \beq\label{a46}
  \begin{array}{c}
  \displaystyle{
 R^\hbar_{12}(z_1,z_2) R^\hbar_{23}(z_2,z_3)
 R^\hbar_{31}(z_3,z_1)+R^\hbar_{13}(z_1,z_3)
 R^\hbar_{32}(z_3,z_2) R^\hbar_{21}(z_2,z_1)=-\wp'(\hbar)\,\,1\otimes 1\otimes 1
 }
 \end{array}
 \eq
via the usage of skew-symmetry (\ref{a41}) and due to
  \beq\label{a47}
  \begin{array}{c}
  \displaystyle{
 \lim\limits_{\hbar=0} \left[\hbar\, R^{\hbox{\tiny{ACF}}}_{12}(\hbar,z_1,z_2|\,u)\right] = 1\otimes 1\,.
 }
 \end{array}
 \eq

\section{IRF-Vertex for ACF $R$-matrix and higher identities}\label{sect3}
\setcounter{equation}{0}

We start with
\noindent \vskip1mm\underline{\em{Proof of Theorem
\ref{theor2}}}:\vskip1mm

Suppose (\ref{a26}) holds true. Then (\ref{a25}) follows from
(\ref{a07}) and (\ref{a10}). Indeed, plugging (\ref{a26}) into
(\ref{a10}) we get for $R^{\hbox{\tiny{F}}}$
  \beq\label{a30}
  \begin{array}{l}
  \displaystyle{
 R^{\hbox{\tiny{F}}}_{12}(\hbar,z_1-z_2|\,u)=
 }
 \\ \ \\
  \displaystyle{
 =g_1^{-1}(z_1,u-\hbar^{(2)})\,
 g_1(z_1+\hbar,u)\,R^{\hbox{\tiny{ACF}}}_{12}(\hbar,z_1,z_2|\,u)\,
 g_2^{-1}(z_2+\hbar,u)\, g_2(z_2,u-\hbar^{(1)})\,.
 }
 \end{array}
 \eq
 On the other hand, it follows from the IRF-Vertex relation (\ref{a07})
 that $R^{\hbox{\tiny{F}}}_{12}(\hbar,z_1-z_2|\,u)$ equals
  \beq\label{a31}
  \begin{array}{l}
  \displaystyle{
g_1^{-1}(z_1,u-\hbar^{(2)})\,
 g_2^{-1}(z_2,u)\,R^{\hbox{\tiny{B}}}_{12}(\hbar,z_1-z_2|\,u)\,
 g_1(z_1,u)\, g_2(z_2,u-\hbar^{(1)})\,.
 }
 \end{array}
 \eq
 Compared together (\ref{a30}) and (\ref{a31}) yield (\ref{a25}).

Let us prove now (\ref{a26}), which can be rewritten in the form
  \beq\label{a32}
  \begin{array}{l}
  \displaystyle{
 g_1(z,u)\,{\bar
 R}^{-1}_{12}(\hbar,z|\,u)=g_1(z+\hbar,u+\hbar^{(2)})\,.

 }
 \end{array}
 \eq
Substituting $g(z,u)$ (\ref{a06}) and ${\bar  R}^{-1}_{12}$
(\ref{a12}) into (\ref{a32}). Then we obtain for its l.h.s.:
  \beq\label{a33}
  \begin{array}{l}
  \displaystyle{
\frac{\vth(\hbar)}{\vth'(0)}\sum\limits_{i,j,k} E_{ij}\otimes
E_{kk}\, g_{ij}(z,u)\phi(\hbar,u_{jk}-\hbar) +
 E_{ij}\otimes
E_{jj}\, g_{ik}(z,u)\phi(z,\hbar-u_{kj})\,.
 }
 \end{array}
 \eq
Using the definition (\ref{a041}), the r.h.s. of (\ref{a32}) equals
  \beq\label{a34}
  \begin{array}{l}
  \displaystyle{
\sum\limits_{i,j,k} E_{ij}\otimes E_{kk}\,
\exp{\left(\hbar\frac{\p}{\p u_k}\right)}\, g_{ij}(z+\hbar,u)\,
\exp{\left(-\hbar\frac{\p}{\p u_k}\right)}\,.
 }
 \end{array}
 \eq
Comparing (\ref{a33}) and (\ref{a34}) taking into account that
$\exp(\hbar\p_{u_k})\left(\sum_m
u_m\right)\exp(-\hbar\p_{u_k})=\hbar+\sum_m u_m$ for any $k$. For
$k\neq j$ type terms equality of (\ref{a33}) and (\ref{a34}) is
equivalent to
  \beq\label{a35}
  \begin{array}{l}
  \displaystyle{
\frac{\vth(\hbar)}{\vth'(0)}\,g_{ij}(z,u)\,\phi(\hbar,u_{jk}-\hbar)=g_{ij}(z,u)\,\frac{\vth(u_{kj})}{\vth(u_{kj}+\hbar)}\,.
 }
 \end{array}
 \eq
It is the definition of the Kronecker function. For $k=j$ the first
term in (\ref{a33}) equals zero ($\phi(\hbar,-\hbar)=0$). The
equality of (\ref{a33}) and (\ref{a34}) takes the form:
  \beq\label{a36}
  \begin{array}{l}
  \displaystyle{
\frac{\vth(\hbar)}{\vth'(0)}\,\sum\limits_k g_{ik}(z,u)\,
\phi(z,-u_{kj}+\hbar)=g_{ij}(z+N\hbar,u) \prod\limits_{m\neq
j}\frac{\vth(u_{mj})}{\vth(u_{mj}-\hbar)}\,.
 }
 \end{array}
 \eq
The latter is the statement of \cite{Hasegawa} about factorization
of the Lax operator for the Ruijsenaars-Schneider model.
$\blacksquare$

%quasi-hopf ???? Jimbo page 3 "in general, the new coproduct is no
%longer coassociative" -- possibly, the requirement of coassoc leads
%to ACF as solution of AYBE

\noindent {\bf Identities for ACF $R$-matrix.} The statement of
Theorem \ref{theor2} allows to derive results of the previous
Section from those for non-dynamical $R$-matrices. Let us show that
(\ref{a23}) follows from the associative Yang-Baxter equation
(\ref{a15}) for non-dynamical $R$-matrices.
To see it we will use that the l.h.s. (i.e. the Baxter-Belavin
$R$-matrix) of (\ref{a25}) depends on difference of spectral
parameters, that is its r.h.s. is independent of the shift
$z_1\rightarrow z_1+c$ and $z_2\rightarrow z_2+c$.

Substitute (\ref{a25}) into
$R^{\hbox{\tiny{B}}}_{12}(\hbar,z_1-z_2)R^{\hbox{\tiny{B}}}_{23}(\eta,z_2-z_3)$
(the l.h.s. of (\ref{a15})) with the constant $c=\eta$ for
$R^{\hbox{\tiny{B}}}_{12}$ and with $c=\hbar$ for
$R^{\hbox{\tiny{B}}}_{23}$:
  \beq\label{a37}
  \begin{array}{c}
  \displaystyle{
R^{\hbox{\tiny{B}}}_{12}(\hbar,z_1-z_2)R^{\hbox{\tiny{B}}}_{23}(\eta,z_2-z_3)=g_1(z_1+\hbar+\eta)g_2(z_2+\eta)g_3(z_3+\hbar)\times
 }
 \\ \ \\
  \displaystyle{
\times
R^{\hbox{\tiny{ACF}}}_{12}(\hbar,z_1+\eta,z_2+\eta)R^{\hbox{\tiny{ACF}}}_{23}(\eta,z_2+\hbar,z_3+\hbar)
 g_1^{-1}(z_1+\eta) g_2^{-1}(z_2+\hbar) g_3^{-1}(z_3+\hbar+\eta)\,.
 }
 \end{array}
 \eq
In the same way substitute (\ref{a25}) into
$R^{\hbox{\tiny{B}}}_{13}(\eta,z_1-z_3)R^{\hbox{\tiny{B}}}_{12}(\hbar-\eta,z_1-z_2)$
(the first term in the r.h.s. of (\ref{a15})) with the constant
$c=\hbar$ for $R^{\hbox{\tiny{B}}}_{13}$ and with $c=\eta$ for
$R^{\hbox{\tiny{B}}}_{12}$:
  \beq\label{a38}
  \begin{array}{c}
  \displaystyle{
R^{\hbox{\tiny{B}}}_{13}(\eta,z_1-z_3)R^{\hbox{\tiny{B}}}_{12}(\hbar-\eta,z_1-z_2)=g_1(z_1+\hbar+\eta)g_2(z_2+\eta)g_3(z_3+\hbar)
 \times
 }
 \\ \ \\
  \displaystyle{
\times
R^{\hbox{\tiny{ACF}}}_{13}(\eta,z_1\!+\!\hbar,z_3\!+\!\hbar)R^{\hbox{\tiny{ACF}}}_{12}(\hbar-\eta,z_1\!+\!\eta,z_2\!+\!\eta)
 g_1^{-1}(z_1+\eta) g_2^{-1}(z_2+\hbar) g_3^{-1}(z_3+\hbar+\eta)\,.
 }
 \end{array}
 \eq
At last substitute (\ref{a25}) into
$R^{\hbox{\tiny{B}}}_{23}(\eta-\hbar,z_2-z_3)R^{\hbox{\tiny{B}}}_{13}(\hbar,z_1-z_3)$
(the second term in the r.h.s. of (\ref{a15})) with the constant
$c=\hbar$ for $R^{\hbox{\tiny{B}}}_{23}$ and with $c=\eta$ for
$R^{\hbox{\tiny{B}}}_{13}$:
  \beq\label{a39}
  \begin{array}{c}
  \displaystyle{
R^{\hbox{\tiny{B}}}_{23}(\eta-\hbar,z_2-z_3)R^{\hbox{\tiny{B}}}_{13}(\hbar,z_1-z_3)=g_1(z_1+\hbar+\eta)g_2(z_2+\eta)g_3(z_3+\hbar)
 \times
 }
 \\ \ \\
  \displaystyle{
\times
R^{\hbox{\tiny{ACF}}}_{23}(\eta-\hbar,z_2\!+\!\hbar,z_3\!+\!\hbar)R^{\hbox{\tiny{ACF}}}_{13}(\hbar,z_1\!+\!\eta,z_3\!+\!\eta)
 g_1^{-1}(z_1\!+\!\eta) g_2^{-1}(z_2\!+\!\hbar) g_3^{-1}(z_3\!+\!\hbar\!+\!\eta)\,.
 }
 \end{array}
 \eq
The products of $g$ matrices in (\ref{a37}), (\ref{a38}),
(\ref{a39}) are the same. Therefore, from (\ref{a37}), (\ref{a38}),
(\ref{a39}) and (\ref{a15}) we get (\ref{a23}).

Similarly, one can get the Yang-Baxter equation (\ref{a08}) from
(\ref{a01}) and more general cubic relation (\ref{a24}) from
(\ref{a19}). In the latter case we get
  \beq\label{a61}
  \begin{array}{c}
  \displaystyle{
\Big[\hbox{eq.
}(\ref{a19})\Big]=g_1(z_1+\eta)g_2(z_2)g_3(z_3-\hbar)\,
 \Big[\hbox{eq. }(\ref{a24})\Big]\, g_1^{-1}(z_1-\hbar) g_2^{-1}(z_2) g_3^{-1}(z_3+\eta)
 }
 \end{array}
 \eq

Finally, in the same way one can verify that the higher order (in
$R$) identities (\ref{a27}) are also valid for ACF $R$-matrix. It
happens because each term of the sum in (\ref{a27}) contains all
distinct indices. Therefore, each matrix $g$ is either cancelled out
or can be removed to the right or to the left. This reasoning shows
that the substitution of (\ref{a25}) into the identity (written for
the Baxter-Belavin $R$-matrix) yields (\ref{a27}) for the ACF
$R$-matrix conjugated by
 $$
g_a(z_a+\hbar)\prod\limits_{c\neq a} g_c(z_c)\,.
 $$

\vskip1mm

\noindent {\bf Modification of bundles as matrix theta function.} It
was shown in \cite{LOZ,LOZ8} that the intertwining matrix
(\ref{a06}) is of the same form in classical mechanics, where it
plays the role of special gauge transformation relating Lax pairs of
Calogero-Moser (Ruijsenaars-Schneider) model and integrable
(relativistic) elliptic top. This approach treats the Lax operator
of an integrable system as section of some bundle over complex curve
(with local coordinate $z$). The gauge transformation (\ref{a06})
changes its characteristic class (e.g. degree of underlying vector
bundle) because it is degenerated at $z=0$:
  \beq\label{a62}
  \begin{array}{c}
  \displaystyle{
\det g(z,u)=c(\tau)\vth(z)\prod\limits_{j>k}\vth(u_j-u_k)\,.
 }
 \end{array}
 \eq
 Such
gauge transformations are called modifications of bundles, and the
gauge equivalence of a set of integrable systems related to
different characteristic classes is called the symplectic Hecke
correspondence. Here we argue that $g(z,u)$ matrix can be considered
as a matrix analogue of theta function (\ref{a902}) in the same way
as $R$-matrix (\ref{a03}) is a matrix analogue of the Kronecker
function (\ref{a907}). See \cite{Pol} and \cite{LOZ9,LOZ10} for
details.

First, notice that by definition (\ref{a06}) $g(z)$ is indeed
$\vth(z)$ in scalar ($N=1$) case. Next, consider (\ref{a25}) written
as follows
  \beq\label{a63}
  \begin{array}{c}
  \displaystyle{
 g_2^{-1}(z_2,u)\,R^{\hbox{\tiny{B}}}_{12}(\hbar,z_1-z_2)=g_1(z_1+\hbar,u)\,
 R^{\hbox{\tiny{ACF}}}_{12}(\hbar,z_1,z_2|\,u)\, g_2^{-1}(z_2+\hbar,u)\,
 g_1^{-1}(z_1,u)\,.
 }
 \end{array}
 \eq
As functions of $z_2$ both parts of (\ref{a63}) have simple poles at
$z_2=0$. Taking residues at $z_2=0$ we get
  \beq\label{a64}
  \begin{array}{c}
  \displaystyle{
 {\breve g}_2(0,u)\,R^{\hbox{\tiny{B}}}_{12}(\hbar,z)=g_1(z+\hbar,u)\,
 \mathcal O_{12}\, g_2^{-1}(\hbar,u)\,
 g_1^{-1}(z,u)\,,
 }
 \end{array}
 \eq
 where
   \beq\label{a65}
  \begin{array}{c}
  \displaystyle{
{\breve g}(0,u)=\res\limits_{z=0}\,g^{-1}(z)
 }
 \end{array}
 \eq
is a matrix analogue of theta constant $1/\vth'(0)$ while $\mathcal
O_{12}$ is the following (degenerated) matrix:
   \beq\label{a66}
  \begin{array}{c}
  \displaystyle{
 \mathcal O_{12}=\res\limits_{z_2=0}\,R^{\hbox{\tiny{ACF}}}_{12}(\hbar,z_1,z_2|\,u)=\sum\limits_{i,j}E_{ii}\otimes
E_{ji}\,.
 }
 \end{array}
 \eq
When $N=1$ all the elements in (\ref{a64}) become scalar, $\mathcal
O_{12}\left.\right|_{N=1}=1$, and we reproduce the definition of the
Kronecker function (\ref{a907}).

\section{Equations for Felder's $R$-matrix}\label{sect4}
\setcounter{equation}{0}

Let us derive equations for the Felder's $R$-matrix (\ref{a05})
which follow from those for the Baxter-Belavin case via the
IRF-Vertex transformation (\ref{a07}). Rewrite (\ref{a07}) as
follows:
 \beq\label{a72}
 \begin{array}{c}
  \displaystyle{
R^{\hbox{\tiny{B}}}_{12}(\hbar,z_1-z_2) =g_2(z_2,u)\,
g_1(z_1,u-\hbar^{(2)})\,R^{\hbox{\tiny{F}}}_{12}(\hbar,z_1-z_2|\,u)\,g_2^{-1}(z_2,u-\hbar^{(1)})\,g_1^{-1}(z_1,u)\,,
 }
 \end{array}
 \eq
 \beq\label{a71}
 \begin{array}{c}
  \displaystyle{
R^{\hbox{\tiny{B}}}_{12}(\hbar,z_1-z_2)=g_1(z_1,u)\,g_2(z_2,u+\hbar^{(1)})\,R^{\hbox{\tiny{F}}}_{12}(\hbar,z_1-z_2|\,u)
\,g_1^{-1}(z_1,u+\hbar^{(2)})\,g_2^{-1}(z_2,u)\,.
 }
 \end{array}
 \eq
The latter follows from (\ref{a07}) and the skew-symmetry
(\ref{a16}) valid for $R^{\hbox{\tiny{B}}}$ and
$R^{\hbox{\tiny{F}}}$. Let us mention that (\ref{a71}) together with
(\ref{a101}) reproduces the relation (\ref{a25}) between the
Baxter-Belavin and ACF $R$-matrices  in the same way as it was shown
in (\ref{a30})-(\ref{a31}) using (\ref{a26}) and (\ref{a07}),
(\ref{a10}).

Plugging (\ref{a72}) or (\ref{a71}) into any term from the
associative Yang-Baxter equation (\ref{a15}) it is easy to see that
one can not cancel out all the multiples of $g$ matrices between a
pair of $R^{\hbox{\tiny{F}}}$. There are only the possibilities to
keep a single multiple of this type between a pair of $R$-matrices.
Consider for example transformation of
$R^{\hbox{\tiny{B}}}_{12}(\hbar,z_1-z_2)R^{\hbox{\tiny{B}}}_{23}(\eta,z_2-z_3)$:
 \beq\label{a73}
 \begin{array}{c}
  \displaystyle{
R^{\hbox{\tiny{B}}}_{12}(\hbar,z_1-z_2)R^{\hbox{\tiny{B}}}_{23}(\eta,z_2-z_3)=
 (P_3^{-\eta}R^{\hbox{\tiny{B}}}_{12}(\hbar,z_1-z_2)P_3^{\eta})R^{\hbox{\tiny{B}}}_{23}(\eta,z_2-z_3)=
 }
 \\ \ \\
  \displaystyle{
=g_3(z_3,u)g_1(z_1,u-\eta^{(3)})g_2(z_2,u+\hbar^{(1)}-\eta^{(3)})
R^{\hbox{\tiny{F}}}_{12}(\hbar,z_1-z_2|\,u-\eta^{(3)})\times
 }
  \\ \ \\
  \displaystyle{
g_1^{-1}(z_1,u+\hbar^{(2)}-\eta^{(3)})R^{\hbox{\tiny{F}}}_{23}(\hbar,z_2-z_3|\,u)
g_3^{-1}(z_3,u-\eta^{(2)})g_2^{-1}(z_2,u)\,.
 }
 \end{array}
 \eq
Here we used (\ref{a71}) for
$R^{\hbox{\tiny{B}}}_{12}(\hbar,z_1-z_2)$ and (\ref{a72}) for
$R^{\hbox{\tiny{B}}}_{23}(\eta,z_2-z_3)$. Another possibility to
keep a single $g$ multiple between $R^{\hbox{\tiny{F}}}$ is
 \beq\label{a74}
 \begin{array}{c}
  \displaystyle{
R^{\hbox{\tiny{B}}}_{12}(\hbar,z_1-z_2)R^{\hbox{\tiny{B}}}_{23}(\eta,z_2-z_3)=
 R^{\hbox{\tiny{B}}}_{12}(\hbar,z_1-z_2)(P_1^{-\hbar}R^{\hbox{\tiny{B}}}_{23}(\eta,z_2-z_3)P_1^{\hbar})=
 }
 \\ \ \\
  \displaystyle{
=g_2(z_2,u)g_1(z_1,u-\hbar^{(2)})
R^{\hbox{\tiny{F}}}_{12}(\hbar,z_1-z_2|\,u)g_3(z_3,u-\hbar^{(1)}+\eta^{(2)})\times
 }
  \\ \ \\
  \displaystyle{
R^{\hbox{\tiny{F}}}_{23}(\hbar,z_2-z_3|\,u-\hbar^{(1)})
g_2^{-1}(z_2,u+\eta^{(3)}-\hbar^{(1)})g_3^{-1}(z_3,u-\hbar^{(1)})\,.
 }
 \end{array}
 \eq
Here we used (\ref{a72}) for
$R^{\hbox{\tiny{B}}}_{12}(\hbar,z_1-z_2)$ and (\ref{a71})
$R^{\hbox{\tiny{B}}}_{23}(\eta,z_2-z_3)$.

Combining application of (\ref{a72}) and (\ref{a71}) to
$R^{\hbox{\tiny{B}}}_{ab}(\hbar,z_a-z_b)$ we can generate identities
for the Felder's $R$-matrix starting from those for
$R^{\hbox{\tiny{B}}}_{ab}(\hbar,z_a-z_b)$, but the $g$ multiples are
not cancelled out from the obtained expressions. Let us however
write down the transformed cubic identity (\ref{a19}). Again we use
that $R^{\hbox{\tiny{B}}}_{12}(\hbar)=P_3^\eta
R^{\hbox{\tiny{B}}}_{12}(\hbar) P_3^{-\eta}$ since it is
non-dynamical:
 \beq\label{a75}
 \begin{array}{c}
  \displaystyle{
R^{\hbox{\tiny{B}}}_{12}(\eta,z_1-z_2)(P_2^{\eta}
R^{\hbox{\tiny{B}}}_{13}(\hbar,z_1-z_3)P_2^{-\eta})R^{\hbox{\tiny{B}}}_{23}(\eta,z_2-z_3)-
 }
 \\ \ \\
  \displaystyle{
-(P_1^{\eta}R^{\hbox{\tiny{B}}}_{23}(\hbar,z_2-z_3)P_1^{-\eta})
R^{\hbox{\tiny{B}}}_{13}(\eta,z_1-z_3)(P_3^{\eta}
R^{\hbox{\tiny{B}}}_{12}(\hbar,z_1-z_2) P_3^{-\eta})=
 }
  \\ \ \\
  \displaystyle{
=(\wp(\eta)-\wp(\hbar))\,R^{\hbox{\tiny{B}}}_{13}(\hbar+\eta,z_1-z_3)
\,.
 }
 \end{array}
 \eq
By applying (\ref{a71})  we obtain:
 \beq\label{a76}
 \begin{array}{l}
  \displaystyle{
R^{\hbox{\tiny{F}}}_{12}(\eta,z_1-z_2|\,u)
g_3(z_3,u+\hbar^{(1)}+\eta^{(2)})
R^{\hbox{\tiny{F}}}_{13}(\hbar,z_1-z_3|\,u+\eta^{(2)})\times
 }
 \\ \ \\
  \displaystyle{
g_1^{-1}(z_1,u+\hbar^{(3)}+\eta^{(2)})
R^{\hbox{\tiny{F}}}_{23}(\eta,z_2-z_3|\,u)\,-
 }
 \\ \ \\
  \displaystyle{
-\,g_3(z_3,u+\hbar^{(2)}+\eta^{(1)})
R^{\hbox{\tiny{F}}}_{23}(\hbar,z_2-z_3|\,u+\eta^{(1)})g_2^{-1}(z_2,u+\hbar^{(3)}+\eta^{(1)})\times
 }
 \\ \ \\
  \displaystyle{
R^{\hbox{\tiny{F}}}_{13}(\eta,z_1-z_3|\,u)
g_2(z_2,u+\hbar^{(1)}+\eta^{(3)})
R^{\hbox{\tiny{F}}}_{12}(\hbar,z_1-z_2|\,u+\eta^{(3)})
g_1^{-1}(z_1,u+\hbar^{(2)}+\eta^{(3)}) =
 }
  \\ \ \\
  \displaystyle{
=(\wp(\eta)-\wp(\hbar))\,g_2^{-1}(z_2,u+\eta^{(1)})g_3(z_3,u+\hbar^{(1)}+\eta^{(1)})R^{\hbox{\tiny{B}}}_{13}(\hbar+\eta,z_1-z_3)\times
 }
  \\ \ \\
  \displaystyle{
g_1^{-1}(z_1,u+\hbar^{(3)}+\eta^{(3)}) g_2(z_2,u+\eta^{(3)})\,.
 }
 \end{array}
 \eq
It can be useful by the following reason. In the case $\hbar=\eta$
its r.h.s. equals zero, and the l.h.s. provides the
Gervais-Neveu-Felder equation (\ref{a04}). It happens because in
this case one can apply the weight zero condition $P_1^\hbar
P_2^\hbar R^{\hbox{\tiny{F}}}_{12}=R^{\hbox{\tiny{F}}}_{12}P_1^\hbar
P_2^\hbar$, and the $g$ multiples cancel out (this is the proof of
the IRF-Vertex transformation). On the other hand, the case
$\hbar=-\eta$ applied to (\ref{a19}) provides (\ref{a20}) which
leads to commutativity of the KZB connections
$[\nabla_a,\nabla_\tau]=0$. Therefore, it would appear reasonable
that (\ref{a76}) leads to commutativity of the dynamical (Felder's)
KZB connections. We will discuss it in our next paper.

%Instead of using IRF-Vertex transformation we can also ...

{\em Remark:} Let us mention that while the associative Yang-Baxter
equation (\ref{a15}) (without any shifts of arguments) is not valid
for the Felder's $R$-matrix, the disappearance error has quite
simple form in the rational case. Consider
$R^\hbar_{ab}=R^{\hbox{\tiny{F}}}_{ab}(\hbar,z_a-z_b|\,u)$
(\ref{a05}) with $\phi(\eta,z)=1/\eta+1/z$, then
  \beq\label{a97}
  \begin{array}{c}
  \displaystyle{
 R^\hbar_{12}
 R^{\eta}_{23}-R^{\eta}_{13}R_{12}^{\hbar-\eta}-R^{\eta-\hbar}_{23}R^\hbar_{13}=
 }
 \\ \ \\
  \displaystyle{
 =\sum_{i \ne j} \frac{1}{(u_{i}-u_j)^2}\, \Big( E_{ij} \otimes E_{jj}
\otimes E_{ji} + E_{ii} \otimes E_{ij} \otimes E_{ji} + E_{ij}
\otimes E_{ji} \otimes E_{ii}-
 }
 \\ \ \\
  \displaystyle{
  - E_{ii} \otimes E_{ii} \otimes E_{jj}
- E_{ii} \otimes E_{jj} \otimes E_{jj} - E_{ii} \otimes E_{jj}
\otimes E_{ii} \Big)\,,
 }
 \end{array}
 \eq
i.e. the r.h.s. is independent of spectral parameters and the Planck
constants.

\section{Appendix}

\subsection{Elliptic functions}
\def\theequation{A.\arabic{equation}}
\setcounter{equation}{0}
 The Riemann theta-functions \cite{Fay} with characteristics on an elliptic curve
 $\Sigma_\tau={\mC^2}/(\mZ\oplus\tau\mZ)$ (with moduli $\tau$,
 Im$\tau>0$) are defined for some integer $N\geq 2$:
 \beq\label{a901}
 \begin{array}{c}
  \displaystyle{
\theta{\left[\begin{array}{c}
a\\
b
\end{array}
\right]}(z|\, \tau ) =\sum_{j\in \mathbb\, Z}
\exp\left(2\pi\imath(j+a)^2\frac\tau2+2\pi\imath(j+a)(z+b)\right)\,,\quad
a\,,b\in\frac{1}{N}\,\mZ\,.
 }
 \end{array}
 \eq
It has the following quasi-periodic properties (i.e. behavior on the
lattice $\mZ\oplus\tau\mZ$):
$$
\theta{\left[\begin{array}{c}
a\\
b
\end{array}
\right]}(z+1|\,\tau )=\exp(2\pi\imath a)\,
\theta{\left[\begin{array}{c}
a\\
b
\end{array}
\right]}(z |\,\tau )\,,
$$
$$
\theta{\left[\begin{array}{c}
a\\
b
\end{array}
\right]}(z+a'\tau|\,\tau ) =\exp\left(-2\pi\imath {a'}^2\frac\tau2
-2\pi\imath a'(z+b)\right) \theta{\left[\begin{array}{c}
a+a'\\
b
\end{array}
\right]}(z|\,\tau )\,.
$$
A shorthand notation for the odd theta-function is used
 \beq\label{a902}
 \begin{array}{c}
  \displaystyle{
\vth(z|\,\tau)\equiv\vth(z)\equiv\theta{\left[\begin{array}{c}
1/2\\
1/2
\end{array}
\right]}(z|\, \tau )\,.
 }
 \end{array}
 \eq
The space of functions (\ref{a901}) is a natural module for the
action of the Heisenberg group \cite{RicheyT}. Its finite
dimensional representation is generated by a pair of matrices
$Q,\Lambda\in\Mat$:
 \beq\label{a903}
 \begin{array}{c}
  \displaystyle{
Q_{kl}=\delta_{kl}\exp(\frac{2\pi
 i}{N}k)\,,\ \ \ \Lambda_{kl}=\delta_{k-l+1=0\,{\hbox{\tiny{mod}}}
 N}\,,\quad Q^N=\Lambda^N=1_{N\times N}\,.
 }
 \end{array}
 \eq
 Then for
 \beq\label{a904}
 \begin{array}{c}
  \displaystyle{
 T_a=T_{a_1 a_2}=\exp\left(\frac{\pi\imath}{N}\,a_1
 a_2\right)Q^{a_1}\Lambda^{a_2}\,,\quad
 a=(a_1,a_2)\in\mZ_N\times\mZ_N
 }
 \end{array}
 \eq
 the following relations hold
  \beq\label{a905}
 \begin{array}{c}
  \displaystyle{
T_\al T_\be=\kappa_{\al,\be} T_{\al+\be}\,,\ \ \
\kappa_{\al,\be}=\exp\left(\frac{\pi \imath}{N}(\be_1
\al_2-\be_2\al_1)\right)\,,
 }
 \end{array}
 \eq
where $\al+\be=(\al_1+\be_1,\al_2+\be_2)$. The permutation operator
takes the form
  \beq\label{a9051}
 \begin{array}{c}
  \displaystyle{
P_{12}=\frac{1}{N}\sum\limits_{\al\in\mZ_N\times\mZ_N} T_\al \otimes
T_{-\al}=\sum\limits_{i,j=1}^N E_{ij}\otimes E_{ji}\,.
 }
 \end{array}
 \eq

\noindent The Kronecker function \cite{Weil} can be defined in terms
of (\ref{a902}):
  \beq\label{a907}
  \begin{array}{l}
  \displaystyle{
 \phi(\eta,z)=\left\{
   \begin{array}{l}
    1/\eta+1/z\quad - \quad \hbox{rational case}\,,
    \\
    \coth(\eta)+\coth(z) \quad - \quad \hbox{trigonometric case}\,,
    \\
    \frac{\vth'(0)\vth(\eta+z)}{\vth(\eta)\vth(z)} \quad - \quad
    \hbox{elliptic
    case}\,.
   \end{array}
 \right.
 }
 \end{array}
 \eq
We also need the first Eisenstein (odd) function and the Weierstrass
(even) $\wp$-function. In rational, trigonometric and elliptic cases
they are given by
  \beq\label{a908}
  \begin{array}{c}
  \displaystyle{
 E_1(z)=\left\{
   \begin{array}{l}
 1/z\,,
\\
   \coth(z)\,,
\\
    \vth'(z)/\vth(z)\,,
   \end{array}
 \right.\hskip12mm  \wp(z)=\left\{
   \begin{array}{l}
 1/z^2\,,
\\
   1/\sinh^2(z)\,,
\\
    -\p_z E_1(z)+\frac{1}{3}\frac{\vth'''(0)}{\vth'(0)}\,.
   \end{array}
 \right.
 }
 \end{array}
 \eq
The properties and identities for (\ref{a907})-(\ref{a908}) can be
found in \cite{Fay}. See also \cite{LOZ10} and the Appendix in
\cite{LOZ11}, where the same notations are used. Here we give only
the most important. It is the Fay trisecant identity
  \beq\label{a909}
  \begin{array}{c}
  \displaystyle{
\phi(\hbar,z)\phi(\eta,w)=\phi(\hbar-\eta,z)\phi(\eta,z+w)+\phi(\eta-\hbar,w)\phi(\hbar,z+w)
 }
 \end{array}
 \eq
 and its degenerations
  \beq\label{a911}
  \begin{array}{c}
  \displaystyle{
 \phi(\eta,z)\phi(\eta,w)=\phi(\eta,z+w)(E_1(\eta)+E_1(z)+E_1(w)-E_1(z+w+\eta))\,,
 }
 \end{array}
 \eq
  \beq\label{a912}
  \begin{array}{c}
  \displaystyle{
 \phi(\hbar,z)\phi(\hbar,-z)=\wp(\hbar)-\wp(z)\,.
 }
 \end{array}
 \eq
The quantum Baxter-Belavin $R$-matrix (\ref{a03}) also uses the set
of (sections of bundle over $\Sigma_\tau$) functions\footnote{The
definition of the Baxter-Belavin $R$-matrix in
\cite{LOZ8}-\cite{LOZ11}, \cite{Z} slightly differs from
(\ref{a03}), (\ref{a910}). The relation is as follows:
$R^\hbar(z):=N R^{N\hbar}(z)$.}
 \beq\label{a910}
 \begin{array}{c}
  \displaystyle{
 \vf_a^\hbar(z)=\exp(2\pi\imath\frac{a_2}{N}\,z)\,\phi(z,\frac{\hbar+a_1+a_2\tau}{N})\,,\quad
 a=(a_1,a_2)\in \mZ_N\times\mZ_N\,.
 }
 \end{array}
 \eq
 %which form a basis in
 %sections of the bundle
 %$\Gamma(\hbox{End}V)$
 %for some vector bundle of degree 1 over $\Sigma_\tau$ (see
 %\cite{LOZ}).

\subsection{Direct proof of associative Yang-Baxter equation\\ for
Arutyunov-Chekhov-Frolov $R$-matrix}

Let us compare the left and the right hand sides of (\ref{a23})
rewriting them with the definition of $R^{\hbox{\tiny{ACF}}}_{13}$
(\ref{a09}). It meant here that in all summands $i \ne j, i \ne k, j
\ne k$.

The cancellation in some components directly follows from definition
(\ref{a907}):
 {\footnotesize
 $$
 \begin{array}{llll}
& E_{ij} \ox E_{ii} \ox E_{ji}: & & 0 = - \phi(z_{13},-u_{ij}) \phi(\eta-\hbar,-u_{ij}) + \phi(\eta-\hbar,-u_{ij}) \phi(z_{13},-u_{ij}), \\
& E_{ij} \ox E_{ii} \ox E_{jj}: & & 0 = \phi(z_1+\hbar+\eta,-u_{ij}) \phi(\eta-\hbar,-u_{ij}) - \phi(\eta-\hbar,-u_{ij}) \phi(z_1+\hbar+\eta,-u_{ij}), \\
& E_{ii} \ox E_{jj} \ox E_{ij}: & & \phi(\hbar,-u_{ij}) \phi(z_3+\hbar,-u_{ij}) = \phi(z_3+\hbar,-u_{ij}) \phi(\hbar,-u_{ij}),\\
& E_{ij} \ox E_{ji} \ox E_{ji}: & -&\phi(z_{12},-u_{ij}) \phi(-z_3-\hbar,-u_{ij}) = - \phi(-z_3-\hbar,-u_{ij}) \phi(z_{12},-u_{ij}), \\
& E_{ii} \ox E_{ji} \ox E_{ij}: & -&\phi(\hbar,-u_{ij}) \phi(z_{32},-u_{ij}) = - \phi(\hbar,-u_{ij}) \phi(z_{32},-u_{ij}), \\
& E_{ij} \ox E_{jj} \ox E_{ij}: & -&\phi(z_1+\hbar+\eta,-u_{ij}) \phi(z_3+\hbar,-u_{ij}) = - \phi(z_3+\hbar,-u_{ij}) \phi(z_1+\hbar+\eta,-u_{ij}), \\
& E_{jj} \ox E_{ij} \ox E_{ji}: & &\phi(z_2+\eta,-u_{ij}) \phi(z_3+\hbar,-u_{ij}) = \phi(z_3+\hbar,-u_{ij}) \phi(z_2+\eta,-u_{ij}), \\
& E_{ii} \ox E_{ji} \ox E_{jj}: & -&\phi(\eta,-u_{ij}) \phi(-z_2-\eta,-u_{ij}) = - \phi(\eta,-u_{ij}) \phi(-z_2-\eta,-u_{ij}), \\
% \end{array}
% $$
% }
% %
% {\footnotesize
% $$
% \begin{array}{llll}
& E_{ij} \ox E_{ij} \ox E_{ji}: & & 0 = \phi(z_{13},-u_{ij}) \phi(z_2+\eta,-u_{ij}) - \phi(z_2+\eta,-u_{ij}) \phi(z_{13},-u_{ij}), \\
& E_{ij} \ox E_{ij} \ox E_{jj}: & & 0 = - \phi(z_1+\hbar+\eta,-u_{ij}) \phi(z_2+\eta,-u_{ij}) + \phi(z_2+\eta,-u_{ij}) \phi(z_1+\hbar+\eta,-u_{ij}), \\
& E_{ij} \ox E_{ji} \ox E_{jj}: & &\phi(z_{12},-u_{ij}) \phi(\eta,-u_{ij}) = \phi(\eta,-u_{ij}) \phi(z_{12},-u_{ij}), \\
& E_{ij} \ox E_{ji} \ox E_{kk}: & &\phi(z_{12},-u_{ij}) \phi(\eta,-u_{ik}) = \phi(\eta,-u_{ik}) \phi(z_{12},-u_{ij}), \\
& E_{jj} \ox E_{ij} \ox E_{kk}: & &\phi(z_2+\eta,-u_{ij}) \phi(\eta,-u_{jk}) = \phi(\eta,-u_{jk}) \phi(z_2+\eta,-u_{ij}), \\
& E_{ij} \ox E_{jk} \ox E_{ki}: & -&\phi(z_1+\hbar+\eta,-u_{ij}) \phi(z_{23},-u_{jk}) = - \phi(z_{23},-u_{jk}) \phi(z_1+\hbar+\eta,-u_{ij}), \\
& E_{ii} \ox E_{jk} \ox E_{kj}: & &\phi(\hbar,-u_{ij}) \phi(z_{23},-u_{jk}) = \phi(z_{23},-u_{jk}) \phi(\hbar,-u_{ij}), \\
% \end{array}
% $$
% }
% %
% {\footnotesize
% $$
% \begin{array}{llll}
& E_{ii} \ox E_{jj} \ox E_{kj}: & &\phi(\hbar,-u_{ij}) \phi(z_3+\hbar,-u_{kj}) = \phi(z_3+\hbar,-u_{kj}) \phi(\hbar,-u_{ij}), \\
& E_{ij} \ox E_{ji} \ox E_{ki}: & &\phi(z_{12},-u_{ij}) \phi(z_3+\hbar,-u_{ki}) = \phi(z_3+\hbar,-u_{ki}) \phi(z_{12},-u_{ij}), \\
& E_{ij} \ox E_{jj} \ox E_{kj}: & -&\phi(z_1+\hbar+\eta,-u_{ij}) \phi(z_3+\hbar,-u_{kj}) = - \phi(z_3+\hbar,-u_{kj}) \phi(z_1+\hbar+\eta,-u_{ij}), \\
& E_{jj} \ox E_{ij} \ox E_{kj}: & &\phi(z_2+\eta,-u_{ij}) \phi(z_3+\hbar,-u_{kj}) = \phi(z_3+\hbar,-u_{kj}) \phi(z_2+\eta,-u_{ij}), \\
& E_{ij} \ox E_{kk} \ox E_{ji}: & & 0 = \phi(z_{13},-u_{ij}) \phi(\hbar-\eta,-u_{jk}) + \phi(\eta-\hbar,-u_{kj}) \phi(z_{13},-u_{ij}), \\
& E_{ij} \ox E_{kk} \ox E_{jj}: & & 0 = - \phi(z_1+\hbar+\eta,-u_{ij}) \phi(\hbar-\eta,-u_{jk}) - \phi(\eta-\hbar,-u_{kj}) \phi(z_1+\hbar+\eta,-u_{ij}), \\
& E_{ij} \ox E_{kj} \ox E_{jj}: & & 0 = - \phi(z_1+\hbar+\eta,-u_{ij}) \phi(z_2+\eta,-u_{kj}) + \phi(z_2+\eta,-u_{kj}) \phi(z_1+\hbar+\eta,-u_{ij}), \\
& E_{ij} \ox E_{kj} \ox E_{jj}: & & 0 = \phi(z_{13},-u_{ij})
\phi(z_2+\eta,-u_{kj}) - \phi(z_2+\eta,-u_{kj})
\phi(z_{13},-u_{ij})\,.
 \end{array}
 $$
 }
Other identities can be proven by direct use of (\ref{a909}):
 {\footnotesize $$ \begin{array}{llll}
& E_{ii} \ox E_{jj} \ox E_{kk}: & &\phi(\hbar, -u_{ij}) \phi(\eta,-u_{jk}) = \\
& & & = \phi(\eta,-u_{ik}) \phi(\hbar-\eta,-u_{ij}) + \phi(\eta-\hbar,-u_{jk}) \phi(\hbar,-u_{ik}), \\
& E_{ij} \ox E_{jj} \ox E_{kk}: & -&\phi(z_1+\hbar+\eta,-u_{ij}) \phi(\eta,-u_{jk}) = \\
& & & = - \phi(\eta,-u_{ik}) \phi(z_1+\hbar,-u_{ij}) + \phi(z_1+\hbar+\eta,-u_{ik}) \phi(z_1+\hbar,-u_{kj}), \\
& E_{ik} \ox E_{kk} \ox E_{ji}: & & 0 = \phi(z_{13},-u_{13}) \phi(z_1+\hbar,-u_{jk}) - \\
& & & - \phi(z_3+\hbar,-u_{ji}) \phi(z_1+\hbar,-u_{ik}) +
\phi(z_3+\hbar,-u_{kj}) \phi(z_{13},-u_{ik}),
 \\
& E_{ij} \ox E_{jk} \ox E_{ki}: & &\phi(z_{12},-u_{ij}) \phi(z_{23},-u_{ik}) = \\
& & & = \phi(z_{13},-u_{ik}) \phi(z_{12},-u_{kj}) + \phi(z_{23},-u_{jk}) \phi(z_{13},-u_{ij}), \\
& E_{jj} \ox E_{ik} \ox E_{kj}: & &\phi(z_2+\eta,-u_{ij}) \phi(z_{23},-u_{jk}) =  \\
& & & = - \phi(z_2+\eta,-u_{ik}) \phi(z_3+\eta,-u_{kj}) + \phi(z_{23},-u_{ik}) \phi(z_3+\eta,-u_{ij}), \\
& E_{ii} \ox E_{jk} \ox E_{kk}: & -&\phi(\hbar,-u_{ij}) \phi(z_2+\hbar+\eta,-u_{jk}) - \phi(z_2+\eta,-u_{ji}) \phi(z_2+\hbar+\eta,-u_{ik}) = \\
& & & = - \phi(z_2+\eta,-u_{jk}) \phi(\hbar,-u_{ik}), \\
& E_{kk} \ox E_{ii} \ox E_{jk}: & & 0 = \phi(\eta-\hbar,-u_{ij}) \phi(z_3+\eta,-u_{jk}) + \\
& & & + \phi(z_3+\hbar,-u_{ji}) \phi(z_3+\eta,-u_{ik}) + \phi(z_3+\hbar,-u_{jk}) \phi(\hbar-\eta,-u_{ki}), \\
& E_{ij} \ox E_{jk} \ox E_{kk}: & -&\phi(z_{12},-u_{ij}) \phi(z_2+\hbar-\eta,-u_{ik}) + \phi(z_1+\hbar+\eta,-u_{ij}) \phi(z_2+\hbar+\eta,-u_{jk}) = \\
& & & = - \phi(z_1+\hbar+\eta,-u_{ik}) \phi(z_{12},-u_{kj})\,.
 \end{array} $$ }
One can rewrite the rest of identities using (\ref{a911}) and prove
them by comparing obtained summands with $E_1$-functions:
 {\footnotesize $$ \begin{array}{llll}
& E_{ii} \ox E_{ij} \ox E_{ji}: & &\displaystyle{\frac{\phi(\hbar,z_2+\eta)}{\phi(\hbar,z_1+\hbar)} }\, \phi(\hbar,z_{12}) \phi(z_{23},-u_{ij}) = \displaystyle{ \frac{\phi(\hbar,z_3+\eta)}{\phi(\hbar,z_1+\eta)} }\, \phi(\hbar,z_{13}) \phi(z_{23},-u_{ij}) - \\
& & & - \phi(z_{13},-u_{ij}) \phi(z_{21},-u_{ij}) + \phi(z_2+\eta,-u_{ij}) \phi(-z_3-\eta,-u_{ij}), \\
& E_{ii} \ox E_{jj} \ox E_{ii}: & -&\phi(\hbar,-u_{ij}) \phi(-\eta,-u_{ij}) = - \phi(z_3+\hbar,-u_{ij}) \phi(-z_3-\eta,-u_{ij}) + \\
& & & + \displaystyle{ \frac{\phi(\eta,z_3+\hbar)}{\phi(\eta,z_1+\hbar)} }\, \phi(\eta,z_{13}) \phi(\hbar-\eta,-u_{ij}) - \displaystyle{\frac{\phi(\hbar,z_3+\eta)}{\phi(\hbar,z_1+\eta)} }\, \phi(\hbar,z_{13}) \phi(\hbar-\eta,-u_{ij}), \\
& E_{ij} \ox E_{jj} \ox E_{ii}: & &\phi(z_1+\hbar+\eta,-u_{ij}) \phi(-\eta,-u_{ij}) = \phi(z_3+\hbar,-u_{ij}) \phi(z_{13},-u_{ij}) - \\
& & &-\displaystyle{\frac{\phi(\eta,z_3+\hbar)}{\phi(\eta,z_1+\hbar)} }\,\phi(\eta,z_{13}) \phi(z_1+\hbar,-u_{ij}), \\
% \end{array}
% $$
% }
% %
% {\footnotesize
% $$
% \begin{array}{llll}
& E_{jj} \ox E_{ij} \ox E_{jj}: & &\phi(-\hbar,-u_{ij}) \phi(z_2+\hbar+\eta,-u_{ij}) + \displaystyle{ \frac{\phi(\eta,z_3+\hbar)}{\phi(\eta,z_2+\hbar)} }\,\phi(\eta,z_{23}) \phi(z_2+\eta,-u_{ij}) = \\
& & & = \phi(z_{23},-u_{ij}) \phi(z_3+\eta,-u_{ij}) - \displaystyle{\frac{\phi(\hbar,z_3+\eta)}{\phi(\hbar,z_1+\eta)} }\, \phi(\hbar,z_{13}) \phi(z_2+\eta,-u_{ij}) + \\
& & & + \displaystyle{\frac{\phi(\eta,z_3+\hbar)}{\phi(\eta,z_1+\hbar)} }\, \phi(\eta,z_{13}) \phi(z_2+\eta,-u_{ij}), \\
& E_{ij} \ox E_{jj} \ox E_{ji}: & &\phi(z_{12},-u_{ij}) \phi(z_{23},-u_{ij}) = \displaystyle{ \frac{\phi(\hbar-\eta,z_2+\eta)}{\phi(\hbar-\eta,z_1+\eta)} }\, \phi(\hbar-\eta,z_{12}) \phi(z_{13},-u_{ij}) + \\
& & & + \phi(z_3+\hbar,-u_{ij}) \phi(-z_1-\hbar,-u_{ij}) +
\displaystyle{
\frac{\phi(\eta-\hbar,z_3+\hbar)}{\phi(\eta-\hbar,z_2+\hbar)} }\,
\phi(\eta-\hbar,z_{23}) \phi(z_{13},-u_{ij}),
 \\
& E_{ij} \ox E_{ji} \ox E_{ij}: & &\phi(z_1+\hbar+\eta,-u_{ij})
\phi(z_{32},-u_{ij}) = \phi(z_{32},-u_{ij})
\phi(z_1+\hbar+\eta,-u_{ij}), \\ \ \\
& E_{jj} \ox E_{ii} \ox E_{ij}: & -&\phi(z_2+\eta,-u_{ij}) \phi(z_{32},-u_{ij}) = - \phi(z_3+\hbar,-u_{ij}) \phi(\eta-\hbar,-u_{ij}) + \\
& & & + \displaystyle{ \frac{\phi(\eta-\hbar,z_3+\hbar)}{\phi(\eta-\hbar,z_2+\hbar)} }\, \phi(\eta-\hbar,z_{23}) \phi(z_3+\hbar,-u_{ij}), \\
& E_{ij} \ox E_{jj} \ox E_{jj}: & -&\phi(z_{12},-u_{ij}) \phi(z_2+\hbar+\eta,-u_{ij}) - \displaystyle{\frac{\phi(\eta,z_3+\hbar)}{\phi(\eta,z_2+\hbar)} }\, \phi(\eta,z_{23}) \phi(z_1+\hbar+\eta,-u_{ij}) = \\
& & & = -\phi(\eta,-u_{ij})\phi(z_1+\hbar,-u_{ij}) - \displaystyle{ \frac{\phi(\hbar-\eta,z_2+\eta)}{\phi(\hbar-\eta,z_1+\eta)} \phi(\hbar-\eta,z_{12}) }\, \phi(z_1+\hbar+\eta,-u_{ij}) - \\
& & & - \displaystyle{ \frac{\phi(\eta-\hbar,z_3+\hbar)}{\phi(\eta-\hbar,z_2+\hbar)} }\, \phi(\eta-\hbar,z_{23}) \phi(z_1+\hbar+\eta,-u_{ij}), \\
& E_{ij} \ox E_{ji} \ox E_{ii}: & -&\phi(z_1+\hbar+\eta,-u_{ij}) \phi(-z_2-\hbar-\eta,-u_{ij}) + \displaystyle{ \frac{\phi(\eta,z_3+\hbar)}{\phi(\eta,z_2+\hbar)} \phi(\eta,z_{23}) }\, \phi(z_{12},-u_{ij}) = \\
& & & = - \phi(z_{32},-u_{ij}) \phi(z_{13},-u_{ij}) + \displaystyle{ \frac{\phi(\eta,z_3+\hbar)}{\phi(\eta,z_1+\hbar)} }\,\phi(\eta,z_{13}) \phi(z_{12},-u_{ij}), \\
 \end{array}
 $$
 }
 {\footnotesize
 $$
 \begin{array}{llll}
& E_{jj} \ox E_{ii} \ox E_{ii}: & &\phi(z_2+\eta,-u_{ij}) \phi(-z_2-\hbar-\eta,-u_{ij}) - \displaystyle{\frac{\phi(\eta,z_3+\hbar)}{\phi(\eta,z_2+\hbar)} \phi(\eta,z_{23}) }\, \phi(-\hbar,-u_{ij}) = \\
& & & = - \displaystyle{\frac{\phi(\eta-\hbar,z_3+\hbar)}{\phi(\eta-\hbar,z_2+\hbar)}\phi(\eta-\hbar,z_{23}) }\, \phi(-\hbar,-u_{ij}) + \phi(-\eta,-u_{ij}) \phi(\eta-\hbar,-u_{ij}), \\
& E_{ii} \ox E_{ii} \ox E_{jj}: & & \displaystyle{ \frac{\phi(\hbar,z_2+\eta)}{\phi(\hbar,z_1+\eta)} \phi(\hbar,z_{12}) }\, \phi(\eta,-u_{ij}) = \phi(\eta-\hbar,-u_{ij}) \phi(\hbar,-u_{ij}) + \\
& & & + \displaystyle{\frac{\phi(\hbar-\eta,z_2+\eta)}{\phi(\hbar-\eta,z_1+\eta)} }\, \phi(\hbar-\eta,z_{12}) \phi(\eta,-u_{ij}) - \phi(z_1+\hbar+\eta,-u_{ij}) \phi(-z_1-\hbar,-u_{ij}), \\
& E_{ii} \ox E_{ij} \ox E_{jj}: & -&\displaystyle{ \frac{\phi(\hbar,z_2+\eta)}{\phi(\hbar,z_1+\eta)}\phi(\hbar,z_{12}) }\, \phi(z_2+\hbar+\eta,-u_{ij}) = \\
& & & = \phi(z_1+\hbar+\eta,-u_{ij}) \phi(z_{21},-u_{ij}) - \phi(z_2+\eta,-u_{ij}) \phi(\hbar,-u_{ij}), \\
& E_{jj} \ox E_{jj} \ox E_{ij}: & & \displaystyle{\frac{\phi(\hbar,z_2+\eta)}{\phi(\hbar,z_1+\eta)}\, \phi(\hbar,z_{12}) }\, \phi(z_3+\hbar,-u_{ij}) = \phi(z_{31},-u_{ij}) \phi(z_1+\hbar,-u_{ij}) + \\
& & & + \displaystyle{\frac{\phi(\hbar-\eta,z_2+\eta)}{\phi(\hbar-\eta,z_1+\eta)} }\, \phi(\hbar-\eta,z_{12}) \phi(z_3+\hbar,-u_{ij}) - \phi(\hbar-\eta,-u_{ij}) \phi(z_3+\eta,-u_{ij}) + \\
& & & +
\displaystyle{\frac{\phi(\hbar,z_3+\eta)}{\phi(\hbar,z_1+\eta)}}\,
\phi(\hbar,z_{13}) \phi(z_3+\hbar,-u_{ij})\,.
 \end{array} $$ }
The last identity
 {\footnotesize $$ \begin{array}{l}
 E_{ii} \ox E_{ii} \ox E_{ii}:\quad\quad  \displaystyle{
\frac{\phi(\hbar,z_2+\eta)
\phi(\eta,z_3+\hbar)}{\phi(\hbar,z_1+\eta) \phi(\eta,z_2+\hbar)} }\,
\phi(\hbar,z_{12}) \phi(\eta,z_{23}) = \\ \ \\
  = \displaystyle{\frac{\phi(\hbar-\eta,z_2+\eta)
\phi(\eta,z_3+\hbar)}{\phi(\hbar-\eta,z_1+\eta)
\phi(\eta,z_1+\hbar)} }\, \phi(\eta,z_{13}) \phi(\hbar-\eta,z_{12})
 + \displaystyle{ \frac{\phi(\eta-\hbar,z_3+\hbar)
\phi(\hbar,z_3+\eta)}{\phi(\eta-\hbar,z_2+\hbar)
\phi(\hbar,z_1+\eta)} }\, \phi(\eta-\hbar,z_{23}) \phi(\hbar,z_{13})
 \end{array} $$ }
follows from (\ref{a909}) and the statement
%
% {\footnotesize
 \beq\label{a99}
\begin{array}{c}
 \displaystyle{ \frac{\phi(\hbar,z_2+\eta) \phi(\eta,z_3+\hbar)}{\phi(\hbar,z_1+\eta) \phi(\eta,z_2+\hbar)} =
 \frac{\phi(\hbar-\eta,z_2+\eta) \phi(\eta,z_3+\hbar)}{\phi(\hbar-\eta,z_1+\eta) \phi(\eta,z_1+\hbar)} }= \\ \ \\ %\notag \\
 = \displaystyle{ \frac{\phi(\eta-\hbar,z_3+\hbar) \phi(\hbar,z_3+\eta)}{\phi(\eta-\hbar,z_2+\hbar) \phi(\hbar,z_1+\eta)} =
 \frac{\vth(z_1+\eta) \vth(z_2+\hbar) \vth(z_3+\eta+\hbar)}{\vth(z_1+\hbar+\eta) \vth(z_2+\eta)
 \vth(z_3+\hbar)}\,. }
\end{array}
 \eq

\paragraph*{Acknowledgments.} We are grateful to J. Avan for bringing
the papers \cite{Avan1,Avan2} to our attention. The work was
supported by RFBR grant 15-02-04175, by joint project 15-51-52031
HHC$_a$, by the Dmitry Zimin's fund "Dynasty" and by the Russian
President grant for support of scientific schools NSh-1500.2014.2.

\begin{small}

\end{small}

\end{document}